\documentclass[a4paper, 11pt, reqno]{amsart}

\usepackage{amsfonts}
\usepackage{amsmath}
\usepackage{amssymb}
\usepackage{graphicx}
\usepackage[usenames]{color}
\usepackage[colorlinks]{hyperref}

\begin{document}

\title[On generalized Mittag-Leffler-type functions ...]
{On generalized Mittag-Leffler-type functions of two variables}

\author[A. Hasanov]{Anvar Hasanov}
\address{
  Anvar Hasanov:
  \endgraf
  V.I.Romanovskiy Institute of Mathematics
   \endgraf
  {\it E-mail address} {\rm anvarhasanov@yahoo.com}
  }

\author[E.T.Karimov]{Erkinjon Karimov}
\address{
  Erkinjon Karimov:
  \endgraf
V.I.Romanovskiy Institute of Mathematics
  \endgraf
  Tashkent, Uzbekistan
  \endgraf
   and
  \endgraf
 Ghent University
  \endgraf
  Ghent, Belgium
  \endgraf
{\it E-mail address} {\rm erkinjon@gmail.com; erkinjon.karimov@ugent.be}
 }

\thanks{Authors are partially supported by Methusalem programme of the Ghent University Special Research Fund (BOF), Grant/Award Number: 01M01021; FWO Odysseus 1, Grant/Award Number: G.0H94.18N}



\begin{abstract}
We aim to study  Mittag-Leffler type functions of two variables ${{D}_{1}}\left( x,y \right),...,{{D}_{5}}\left( x,y \right)$ by analogy with the Appell hypergeometric functions of two variables. Moreover, we targeted functions ${{E}_{1}}\left( x,y \right),$ $...,{{E}_{10}}\left( x,y \right)$  as limiting cases of the functions ${{D}_{1}}\left( x,y \right),$ $...,{{D}_{5}}\left( x,y \right)$ and studied certain properties, as well. Following Horn's method, we determine all possible cases of the convergence region of the function ${{D}_{1}}\left( x,y \right).$ Further, for a generalized hypergeometric function, ${{D}_{1}}\left( x,y \right)$ (two variable Mittag-Leffler-type function) integral representations of the Euler type have been proved. One-dimensional and two-dimensional Laplace transforms of the function are also defined. We have constructed a system of partial differential equations which is linked with the function ${{D}_{1}}\left( x,y \right)$.
\end{abstract}

\maketitle

\section{Introduction}

\subsection{Single variable Mittag-Leffler functions and their generalizations.}

The Mittag-Leffler function became popular through its applications \cite{1}. Namely, it appears as a solution of fractional differential equations and integral equations of fractional order. For instance, the Mittag-Leffler function of many variables arises in solving some boundary value problems involving fractional Volterra type integrodifferential equations \cite{3}, initial-boundary value problems for a generalized polynomial diffusion equation with fractional time \cite{4}, and also initial-boundary value problems for time-fractional diffusion equations with positive constant coefficients  \cite{5}. Moreover, the Mittag-Leffler function plays an important role in various fields of applied mathematics and engineering sciences, such as chemistry, biology, statistics, thermodynamics, mechanics, quantum physics, computer science, and signal processing\cite{2}. 

The definition of the classical Mittag-Leffler function is the following \cite{1}:
\begin{equation*}
{E_\alpha }\left( z \right) = \sum\limits_{n = 0}^\infty  {\frac{{{z^n}}}{{\Gamma \left( {\alpha n + 1} \right)}}} ,\,\,\left( {\alpha  > 0,\,z \in \mathbb{C}} \right).
\end{equation*}
In \cite{6}, a function with two parameters was introduced:
\begin{equation*}
{E_{\alpha ,\beta }}\left( z \right) = \sum\limits_{n = 0}^\infty  {} \frac{{{z^n}}}{{\Gamma \left( {\alpha n + \beta } \right)}},\,\,\left( {\alpha ,\beta  \in \mathbb{C} ,\,\,{\mathop{\rm Re}\nolimits} \left( \alpha  \right) > 0,\,\,{\mathop{\rm Re}\nolimits} \left( \beta  \right) > 0} \right).
\end{equation*}
Prabhakar \cite{7} considered a function $E_{\alpha ,\beta }^\gamma \left( z \right)$  with three parameters
\begin{equation*}
\begin{array}{l}
\displaystyle{E_{\alpha ,\beta }^\gamma \left( z \right) = \sum\limits_{n = 0}^\infty  {} \frac{{{{\left( \gamma  \right)}_n}{z^n}}}{{\Gamma \left( {\alpha n + \beta } \right)n!}},\,\,}\\
\left( {\alpha ,\beta ,\gamma  \in \mathbb{C} ,\,\,{\mathop{\rm Re}\nolimits} \left( \alpha  \right) > 0,\,\,{\mathop{\rm Re}\nolimits} \left( \beta  \right) > 0,\,\,{\mathop{\rm Re}\nolimits} \left( \gamma  \right) > 0} \right).
\end{array}
\end{equation*}
The article \cite{8} considers the function
\begin{equation*}
\begin{array}{l}
E_{\alpha ,\beta }^{\gamma ,q}\left( z \right) = \sum\limits_{n = 0}^\infty  {} \displaystyle \frac{{{{\left( \gamma  \right)}_{nq}}{z^n}}}{{\Gamma \left( {\alpha n + \beta } \right)n!}},\,\,\\
\left( {\alpha ,\beta ,\gamma  \in \mathbb{C},\,\,{\mathop{\rm Re}\nolimits} \left( \alpha  \right) > 0,\,\,{\mathop{\rm Re}\nolimits} \left( \beta  \right) > 0,\,\,{\mathop{\rm Re}\nolimits} \left( \gamma  \right) > 0,\,q \in \left( {0,1} \right)} \right).
\end{array}
\end{equation*}
In \cite{9}, \cite{10}, the properties of the following functions were studied
\begin{equation*}
\begin{array}{l}
E_{\alpha ,\beta }^{\gamma ,\delta }\left( z \right) = \sum\limits_{n = 0}^\infty  {} \displaystyle  \frac{{{{\left( \gamma  \right)}_n}{z^n}}}{{\Gamma \left( {\alpha n + \beta } \right){{\left( \delta  \right)}_n}}},\,\,\\
\left( {\alpha ,\beta ,\gamma  \in \mathbb{C},\,\,{\mathop{\rm Re}\nolimits} \left( \alpha  \right) > 0,\,\,{\mathop{\rm Re}\nolimits} \left( \beta  \right) > 0,\,\,{\mathop{\rm Re}\nolimits} \left( \gamma  \right) > 0,\,\delta  > 0} \right),
\end{array}
\end{equation*}
\begin{equation*}
\begin{array}{l}
E_{\alpha ,\beta ,p}^{\gamma ,\delta ,q}\left( z \right) = \sum\limits_{n = 0}^\infty  \displaystyle \frac{{{{\left( \gamma  \right)}_{nq}}{z^n}}}{{\Gamma \left( {\alpha n + \beta } \right){{\left( \delta  \right)}_{np}}}},\,\,\alpha ,\beta ,\gamma  \in \mathbb{C},\,\,{\mathop{\rm Re}\nolimits} \left( \alpha  \right) > 0,\,\,{\mathop{\rm Re}\nolimits} \left( \beta  \right) > 0,\,\,\\
{\mathop{\rm Re}\nolimits} \left( \gamma  \right) > 0,\,\delta  > 0,\,\left( {p,q} \right) > 0,q \le {\mathop{\rm Re}\nolimits} \left( \alpha  \right) + p.
\end{array}
\end{equation*}
The following generalized Mittag-Leffler functions were introduced and studied in the articles \cite{11}, \cite{12}
\begin{equation*}
\begin{array}{l}
{E_{\gamma ,K}}\left[ {{{\left( {{\alpha _j},{\beta _j}} \right)}_{1,m}};z} \right] =\displaystyle \sum\limits_{r = 0}^\infty  \frac{{{{\left( \gamma  \right)}_{rK}}{z^r}}}{{\prod\nolimits_{j = 1}^m {} \Gamma \left( {{\alpha _j}r + {\beta _j}} \right)r!}},\,\,\\
\left( {z,\gamma ,{\alpha _j},{\beta _j} \in \mathbb{C},\,\,\sum\nolimits_{j = 1}^m {} {\mathop{\rm Re}\nolimits} \left( {{\alpha _j}} \right) > {\mathop{\rm Re}\nolimits} \left( K \right) - 1,\,\,j = 1,...,{\mathop{\rm Re}\nolimits} \left( K \right) > 0} \right),
\end{array}
\end{equation*}
\begin{equation*}
\begin{array}{l}
E_{\left( {{\rho _1},...,{\rho _m}} \right),\lambda }^{\left( {{\gamma _1},...,{\gamma _m}} \right)}\left( {{z_1},...,{z_m}} \right) = \displaystyle \sum\limits_{{k_1},...,{k_m} = 0}^\infty  {} \frac{{{{\left( {{\gamma _1}} \right)}_{{k_1}}}{{\left( {{\gamma _2}} \right)}_{{k_2}}} \cdot  \cdot  \cdot {{\left( {{\gamma _k}} \right)}_{{k_m}}}z_1^{{k_1}}z_2^{{k_2}} \cdot  \cdot  \cdot z_m^{{k_m}}}}{{\Gamma \left[ {\lambda  + \sum\nolimits_{j = 1}^m {} {\rho _i}{k_i}} \right]{k_1}!{k_2}! \cdot  \cdot  \cdot {k_m}!}},\,\,\\
\left( {\lambda ,{\gamma _j},{\rho _j},{z_j} \in \mathbb{C},\,{\mathop{\rm Re}\nolimits} \left( {{\rho _j}} \right) > 0,\,j = 1,...,m} \right).
\end{array}
\end{equation*}

\subsection{Comparison analysis.}

In this subsection, we will conduct a short survey of known multivariable, precisely, two-variable Mittag-Leffler type functions. Let us start with a multivariate analog of the Mittag-Leffler function  $E_{\alpha}\left(z\right)$ which was proposed by Luchko Y. and Gorenflo R. \cite{13}. The authors used an operational method to solve a boundary value problem for linear fractional differential equations with constant coefficients. The following function expresses the solution to the studied boundary value problem:
\begin{equation}\label{eq 1.9}
\begin{array}{l}
{E_{\left( {{\alpha _1},...,{\alpha _m}} \right),\beta }}\left( {{z_1},...,{z_m}} \right) = \\
 \sum\limits_{k = 0}^\infty  {} \sum\limits_{
 \begin{array}{l}
 \scriptstyle{l_1} + ... + {l_m} = k\hfill\\
\scriptstyle{l_1} \ge 0,...{l_m} \ge 0\hfill
\end{array}
}^\infty  {} \dfrac{{k!}}{{{l_1}! \times ... \times {l_m}!}}  \dfrac{{\prod\nolimits_{i = 1}^m {z_i^{{l_i}}} }}{{\Gamma \left( {\beta  + \displaystyle \sum\nolimits_{j = 1}^m {} {\alpha _i}{l_i}} \right)}}.
\end{array}
\end{equation}
Certain properties, namely, some recurrence formulas and the important estimate for this function are studied in \cite{5}. In recent work \cite{14}, authors studied some properties, including integral representations, operational relations, differential and pure recurrence relations, Euler, Mellin, and Whittaker transforms of the function (\ref{eq 1.9}) in a particular case, i.e. $m=2$.  

Srivastava H.M., Daoust Martha C. \cite{15}, \cite{16} studied the domain of convergence and Euler-type integral representations for generalized Kampe de Feriet's functions
\begin{equation*}
\begin{array}{l}
S_{C:D;D'}^{A:B;B'}\left( {\begin{array}{*{20}{c}}
x\\
y
\end{array}} \right) = S_{C:D;D'}^{A:B;B'}\left( {\begin{array}{*{20}{c}}
{\left[ {\left( a \right):\theta ,\phi } \right]:\left[ {\left( b \right):\psi } \right];\left[ {\left( {b'} \right):\psi '} \right];}\\
{\left[ {\left( c \right):\delta ,\varepsilon } \right]:\left[ {\left( d \right):\eta } \right];\left[ {\left( {d'} \right):\eta '} \right];}
\end{array}x,y} \right)\\
= \displaystyle \sum\limits_{m,n = 0}^\infty  {\frac{{\prod\limits_{j = 1}^A {} \Gamma \left[ {{a_j} + m{\theta _j} + n{\phi _j}} \right]\prod\limits_{j = 1}^B {} \Gamma \left[ {{b_j} + m{\psi _j}} \right]\prod\limits_{j = 1}^{B'} {} \Gamma \left[ {b{'_j} + n\psi {'_j}} \right]}}{{\prod\limits_{j = 1}^C {} \Gamma \left[ {{c_j} + m{\delta _j} + n{\varepsilon _j}} \right]\prod\limits_{j = 1}^D {} \Gamma \left[ {{d_j} + m{\eta _j}} \right]\prod\limits_{j = 1}^{D'} {} \Gamma \left[ {d{'_j} + n\eta {'_j}} \right]}}\frac{{{x^m}}}{{m!}}\frac{{{y^n}}}{{n!}}},\\
{a_j},{b_j},b{'_j},{c_j},{d_j},d{'_j} \in \mathbb{C} ,\,\,\,{\theta _j},{\phi _j},{\psi _j},\psi {'_j},\,{\delta _j},{\varepsilon _j},{\eta _j},\eta {'_j} \in \mathbb{R}{^ + }.
\end{array}
\end{equation*}

 In \cite{17} the Mittag-Leffler type function  ${E_1}$ of two variables 
 \begin{equation*}
\begin{array}{l}
{E_1}\left( {\begin{array}{*{20}{c}}
{{\gamma _1},{\alpha _1};{\gamma _2},{\beta _1};}\\
{{\delta _1},{\alpha _2},{\beta _2};{\delta _2},{\alpha _3};{\delta _3},{\beta _3};}
\end{array}\left| {\begin{array}{*{20}{c}}
x\\
y
\end{array}} \right.} \right) = \\
\displaystyle \sum\limits_{m,n = 0}^\infty  {\frac{{{{\left( {{\gamma _1}} \right)}_{{\alpha _1}m}}{{\left( {{\gamma _2}} \right)}_{{\beta _1}n}}}}{{\Gamma \left( {{\delta _1} + {\alpha _2}m + {\beta _2}n} \right)}}\frac{{{x^m}}}{{\Gamma \left( {{\delta _2} + {\alpha _3}m} \right)}}\frac{{{y^n}}}{{\Gamma \left( {{\delta _3} + {\beta _3}n} \right)}}} ,\\
{\gamma _1},{\gamma _2},{\delta _1},{\delta _2},{\delta _3},x,y \in \mathbb{C} ,\,\,\,\min \left\{ {{\alpha _1},{\alpha _2},{\alpha _3},{\beta _1},{\beta _2},{\beta _3}} \right\} > 0
\end{array}
\end{equation*}
 was introduced and studied, which in a particular case includes several Mittag-Leffler-type functions of a single variable. The region of convergence determines all possible cases. The system of hypergeometric equations is determined, which satisfies the function ${E_1}$, Euler type integral representations, the Mellin-Barnes contour integral, and the Laplace integral transformation was presented. In that paper, another two-variable Mittag-Leffler type function was also introduced, but not studied:
\begin{equation}\label{eq 1.12}
\begin{array}{l}
{E_2}\left( {\begin{array}{*{20}{c}}
{{\gamma _1},{\alpha _1},{\beta _1};{\gamma _2},{\alpha _2};}\\
{{\delta _1},{\alpha _3},{\beta _2};{\delta _2},{\alpha _4};{\delta _3},{\beta _3};}
\end{array}\left| {\begin{array}{*{20}{c}}
x\\
y
\end{array}} \right.} \right) =\\
 \displaystyle \sum\limits_{m,n = 0}^\infty  {\frac{{{{\left( {{\gamma _1}} \right)}_{{\alpha _1}m + {\beta _1}n}}{{\left( {{\gamma _2}} \right)}_{{\alpha _2}m}}}}{{\Gamma \left( {{\delta _1} + {\alpha _3}m + {\beta _2}n} \right)}}\frac{{{x^m}}}{{\Gamma \left( {{\delta _2} + {\alpha _4}m} \right)}}\frac{{{y^n}}}{{\Gamma \left( {{\delta _3} + {\beta _3}n} \right)}}} ,\\
{\gamma _1},{\gamma _2},{\delta _1},{\delta _2},{\delta _3},x,y \in \mathbb{C},\,\,\,\min \left\{ {{\alpha _1},{\alpha _2},{\alpha _3},{\alpha _4},{\beta _1},{\beta _2},{\beta _3}} \right\} > 0.
\end{array}
\end{equation}
Certain properties of this function were studied in \cite{18},\cite{19}, \cite{20}. Namely, Euler-type integral representations and estimations were obtained. The following statement holds:

{\bf Lemma 1.}\cite{18}
Let $\Re(\delta_1)>\Re(\gamma_1)>0$. If $\alpha_3=\alpha_1$ and $\beta_2=\beta_1$, then the following integral representation holds true:
\begin{equation*}
	\begin{array}{l}
		\displaystyle{E_2\left(
			\begin{matrix}
				\gamma_1,\alpha_1,\beta_1;\gamma_2,\alpha_2\,&\,|x\\
				\delta_1,\alpha_3,\beta_2;\delta_2, \alpha_4;\delta_3,\beta_3\,&\,|y
			\end{matrix}
			\right)=}\\
		 \displaystyle{=\frac{1}{\Gamma(\gamma_1)\Gamma(\delta_1-\gamma_1)}\int\limits_0^1\xi^{\gamma_1-1}(1-\xi)^{\delta_1-\gamma_1-1}E_{\alpha_4,\delta_2}^{\gamma_2,\alpha_2}\left(x\xi^{\alpha_1}\right)E_{\beta_3,\delta_3}\left(y\xi^{\beta_1}\right)d\xi.  }\end{array}
\end{equation*}
Here $E_{\beta_3,\delta_3}(z)$ is two-parameter Mittag-Leffler function and
\[
E_{\alpha_4,\delta_2}^{\gamma_2,\alpha_2}(z)=\sum\limits_{m=0}^\infty\frac{(\gamma_2)_{\alpha_2 m}z^m}{\Gamma(\alpha_4m+\delta_2)}.
\]
If $\lambda<0,\,\delta\le 0$, then the following estimate holds true \cite{18}:
\[
\left|
E_2\left(
\begin{matrix}
	\gamma,\gamma, 1; 1, 0\,&\,|\lambda x^\beta\\
	\beta+k+1,\beta,\alpha;\gamma, \gamma; 1, 1\,&\,|\delta x^\alpha
\end{matrix}
\right)
\right|\le \]
\[\le\frac{1}{\Gamma(\gamma)\Gamma(\beta+k+1-\gamma)}\int\limits_0^1\xi^{\gamma-1}(1-\xi)^{\beta+k-\gamma}\frac{C_1C_2}{1+|\lambda x^\beta\xi^\gamma|}d\xi\le
\]
\[
\le \frac{C_1C_2}{\Gamma(\gamma)\Gamma(\beta+k+1-\gamma)}\int\limits_0^1\xi^{\gamma-1}(1-\xi)^{\beta+k-\gamma}d\xi=\frac{C_1C_2}{\Gamma(\beta+k+1)}=C,
\]
where $C$ is any positive real number.

In \cite{19}, another Euler-type integral representation was presented:

{\bf Lemma 2. } If $\left( {{\alpha }_{3}}>{{\alpha }_{4}},{{\alpha }_{3}}>0,{{\beta }_{2}}>{{\beta }_{3}},{{\beta }_{2}}>0 \right)$, the the following Euler-type integral representations holds true:
\[{{E}_{2}}\left(  
\begin{array}{l}
	\left.{{\gamma }_{1}},{{\alpha }_{1}},{{\beta }_{1}},{{\gamma }_{2}},{{\alpha }_{2}} \right| x\\ 
	\left.{{\delta }_{1}},{{\alpha }_{3}},{{\beta }_{2}},{{\delta }_{2}},{{\alpha }_{4}},{{\delta }_{3}},{{\beta }_{3}}\right| y \\ 
\end{array}\right) =\frac{1}{\Gamma({{\gamma }_{1}})}\cdot \frac{1}{\Gamma({{\gamma }_{2}})}\times
\]
\[
\int\limits_{0}^{1}{\int\limits_{0}^{\infty }{{{\xi }^{\frac{{{\delta }_{1}}}{2}-1}}\cdot {{\eta }^{{{\gamma }_{1}}-1}}\cdot }}{{(1-\xi )}^{\frac{{{\delta }_{1}}}{2}-1}}\cdot {{e}^{-\eta }}e_{{{\alpha }_{3}},-{{\alpha }_{4}}}^{\frac{{{\delta }_{1}}}{2},{{\delta }_{2}}}(x{{\eta }^{{{\alpha }_{1}}}}{{\xi }^{{{\alpha }_{3}}}})\cdot e_{{{\beta }_{2}},-{{\beta }_{3}}}^{\frac{{{\delta }_{1}}}{2},{{\delta }_{3}}}(y{{\eta }^{{{\beta }_{1}}}}{{(1-\xi )}^{{{\beta }_{2}}}})d\xi d\eta, \]  
where 
$e_{\alpha ,\beta }^{\mu ,\delta }(z)=\sum\limits_{k=0}^{\infty }{\frac{{{z}^{n}}}{\Gamma(\alpha k+\mu )\Gamma(\delta -\beta k)}},\alpha >\beta ,\alpha >0$
is the Wright-type function \cite{21}. Then using this representation another estimation of $E_2$ has been presented:

 \textbf{Lemma 3.}\cite{20} If all condition of the Lemma 2 and $$\frac{\delta_{1}}{2}>\alpha_{3}\theta,\frac{\delta_{1}}{2}>\beta_{2}\theta,0\le \theta<\frac{1}{2},\gamma_{1}>\theta(\alpha_{1}+\beta_{1})$$
are valid, then the following holds:
\[
\begin{array}{l}
\left| E_{2} \left(\left. \begin{array}{l} {\gamma_{1} ,\alpha_{1} ,\beta_{1},1,0} \\ {\delta_{1},\alpha_{3} ,\beta_{2},\delta_{2} ,\alpha_{4} ,\delta_{3},\beta_{3}} \end{array}\right|\begin{array}{c} {x} \\ y \end{array}\right) \right|\le\\
\dfrac{{\widetilde{C}}}{(-x)^\theta \cdot(-y^{\theta})}\cdot \dfrac{\Gamma(\frac{\delta_{1}}{2}-\alpha_{3}\theta)\Gamma(\frac{\delta_{1}}{2}-\beta_{2}\theta)\Gamma(\gamma_1-\theta(\alpha_{1}+\beta_{1}))}{\Gamma(\delta_{1}-\theta(\alpha_3+\beta_2))}.
\end{array}
\]

In \cite{22}-\cite{24} the following bivariate Mittag-Leffler functions were introduced:
\[
E_{B,C}^{(\alpha,\beta_1,\beta_2,\gamma_2,\gamma_1)}(x,y)=\sum\limits_{m,n=0}^\infty \frac{(B)_n(C)_m}{\Gamma(\alpha n+\beta_1 m+\gamma_1)\Gamma(\beta_2m+\gamma_2)}\frac{x^n y^m}{n!m!},\,\,\Re(\alpha),\,\Re(\beta)>0,
\]
\[
E_{\gamma_1,\gamma_2,\beta_1,\beta_2}^{(A)}(x,y)=\sum\limits_{m,n=0}^\infty \frac{(A)_{n+m}}{\Gamma(\gamma_1 +\beta_1 n)\Gamma(\gamma_2+\beta_2 m)}\frac{x^n y^m}{n!m!},\,\,\Re(\beta_1),\,\Re(\beta_2)>0,
\]
\[
E_{\alpha,\beta, \gamma}^{A}(x,y)=\sum\limits_{m,n=0}^\infty \frac{(A)_{n+m}}{\Gamma(\alpha n +\beta m+\gamma)}\frac{x^n y^m}{n!m!},\,\,\Re(\alpha),\,\Re(\beta)>0.
\]

In \cite{25}, authors used introduced general bivariate Mittag-Leffler
functions to define fractional integral operators (which have a semigroup property) and corresponding fractional derivative operators (which act as left inverses and analytic continuations). They also showed how these functions and operators arise naturally from some fractional partial integrodifferential equations of the Riemann–Liouville type.

\subsection{Applications of two-variable Mittag-Leffler-type functions.} 

As we mentioned before, two variable Mittag-Leffler-type functions appear in the solution for multi-term fractional order differential equations \cite{13} and mixed equations \cite{26}. 

The paper \cite{26} considered mixed wave-diffusion equation
\begin{equation}\label{eq 1.13}
\begin{array}{l}
f\left( x \right) = \left\{ {\begin{array}{*{20}{c}}
{{}_CD_{0t}^\alpha u\left( {x,t} \right) - {u_{xx}}\left( {x,t} \right),\,\,\,t > 0,}\\
{{}_CD_{t0}^\beta u\left( {x,t} \right) - {u_{xx}}\left( {x,t} \right),\,\,\,t < 0,}
\end{array}} \right.
\end{array}
\end{equation}
in a domain $ \Omega = \{(x,t):0 < x < 1,- p<t <q\}$, where $\alpha,\beta,p,q\in\mathbb{R}{^+},\,0 < \alpha  < 1,\,\,1 < \beta < 2.$
For equation (\ref{eq 1.13}), the boundary value problem is considered and the solution to this problem is expressed by the function ${E_1}$.

A boundary problem in a domain $\Omega  = \left\{ {\left( {x,t} \right):0 < x < 1,\,\,0 < t < T} \right\}$  for the diffusion equation with a fractional time derivative is considered \cite{18}
\begin{equation*}
{}^{PC}D_{0t}^{\alpha ,\beta ,\gamma ,\delta }u\left( {t,x} \right) - {u_{xx}}\left( {t,x} \right) = f\left( {t,x} \right),\,\,\alpha ,\beta ,\gamma ,\delta  \in \mathbb{C} ,\,\,\,{\mathop{\rm Re}\nolimits} \alpha  > 0,
\end{equation*}
where\\
\begin{equation*}
\begin{array}{l}
{}^{PC}D_{0t}^{\alpha ,\beta ,\gamma ,\delta }y\left( t \right) = \displaystyle {}^PI_{0t}^{\alpha ,m - \beta , - \gamma ,\delta }\frac{{{d^m}}}{{d{t^m}}}y\left( t \right),\,\,m - 1 < {\mathop{\rm Re}\nolimits} \beta  < m,\,\,m \in \mathbb{N},\\
{}^PI_{0t}^{\alpha ,\beta ,\gamma ,\delta }y\left( t \right) = \displaystyle \int\limits_0^t {{{\left( {t - \xi } \right)}^{\beta  - 1}}} E_{\alpha ,\beta }^\gamma \left[ {\delta {{\left( {t - \xi } \right)}^\alpha }} \right]y\left( \xi  \right)d\xi.
\end{array}
\end{equation*}
The solution to the problem is expressed in terms of function ${E_2}$. Solutions of direct and inverse problems for the following mixed equation
\begin{equation*} 
	\begin{array}{l}
	\dfrac{1-sign(t-a)}{2} \cdot {}^{PC} D_{ot}^{\alpha ,\beta _{1} ,\gamma ,\delta } u(t,x)+\\
	+\dfrac{1+sign(t-a)}{2} \cdot {}^{PC} D_{at}^{\alpha ,\beta _{2} ,\gamma ,\delta } u(t,x)-u_{xx} (t,x)=f(t,x) 
	\end{array}
\end{equation*}
contain two variable Mittag-Leffler function $E_2$, as well \cite{19}. The solution of an initial-boundary problem for mixed time-fractional equation with space variable coefficients involving the Caputo-Prabhakar fractional derivative is also presented via the function $E_2$ \cite{20}. Also, the solution to the boundary-value problem for the fractional order telegraph equation 
\[
\frac{\partial}{\partial x}{}^{PC}D_{0t}^{\alpha, \beta, \gamma, \delta} u(t,x)+\lambda u(t,x)=f(t,x)
\]
is represented via the variable Mittag-Leffler-type function \cite{27}.  

Another possible application of two variable Mittag-Leffler-type functions is linked to the consideration of integral and differential operators with such functions in the kernel. As shown in \cite{25}, they will be connected with partial integrodifferential equations of the Riemann-Liouvillee type. 

We note that special functions are closely related to fractional calculi (see \cite{28}-\cite{33}), as well as to generalized fractional calculi (see, for example, \cite{34} - \cite{37}). Special functions can be represented as fractional order integration or differentiation operators of some basic elementary special functions. Relations of this kind also provide some alternative definitions for special functions. An example of such unified approaches to special functions can be seen in Kiryakova (\cite{33}, ch. 4). Many recent works on special functions and their application in solving problems from control theory, mechanics, physics, engineering, economics, etc. can be found in the specialized journal "Fractional Calculus and Applied Analysis" (ISSN 1311-0454), available at the website http://www.math.bas.bg/~fcaa. 

We also note interesting targets in special functions such as Le Roy functions. Authors consider multi-index Le Roy functions and introduce different Le Roy functions of Mittag-Leffler-type functions, for instance, Prabhakar type (see \cite{38}-\cite{41}). We think that two variable Mittag-Leffler-type functions might be targeted in a such way.

\section{Definitions}

Having carefully studied the definitions of generalized hypergeometric functions and the Mittag-Leffler type functions, we understand the similarity of these functions with the Horn functions \cite{34}. Given this situation, we define the following functions. Note that the parameters introduced in the functions satisfy the conditions ${\gamma _i},{\delta _i},x,y \in \mathbb{C},$ and ${\alpha _i},{\beta _i} \in \mathbb{R},\,\,\,\,\min \left\{ {{\alpha _i},{\beta _i}} \right\} > 0$:
\begin{equation}\label{eq 2.1}
\begin{array}{l}
{D_1}\left( {\begin{array}{*{20}{c}}
{{\gamma _1},{\alpha _1},{\beta _1};{\gamma _2},{\alpha _2};{\gamma _3},{\beta _2};}\\
{{\delta _1},{\alpha _3},{\beta _3};{\delta _2},{\alpha _4};{\delta _3},{\beta _4};}
\end{array}\left| {\begin{array}{*{20}{c}}
x\\
y
\end{array}} \right.} \right) \\
= \displaystyle \sum\limits_{m,n = 0}^\infty  {\frac{{{{\left( {{\gamma _1}} \right)}_{{\alpha _1}m + {\beta _1}n}}{{\left( {{\gamma _2}} \right)}_{{\alpha _2}m}}{{\left( {{\gamma _3}} \right)}_{{\beta _2}n}}}}{{\Gamma \left( {{\delta _1} + {\alpha _3}m + {\beta _3}n} \right)}}\frac{{{x^m}}}{{\Gamma \left( {{\delta _2} + {\alpha _4}m} \right)}}\frac{{{y^n}}}{{\Gamma \left( {{\delta _3} + {\beta _4}n} \right)}}} ,
\end{array}
\end{equation}
\begin{equation}\label{eq 2.2}
\begin{array}{l}
{D_2}\left( {\begin{array}{*{20}{c}}
{{\gamma _1},{\alpha _1},{\beta _1};{\gamma _2},{\alpha _2};{\gamma _3},{\beta _2};}\\
{{\delta _1},{\alpha _3};{\delta _2},{\beta _3};{\delta _3},{\alpha _4};{\delta _4},{\beta _4};}
\end{array}\left| {\begin{array}{*{20}{c}}
x\\
y
\end{array}} \right.} \right) \\
=\displaystyle \sum\limits_{m,n = 0}^\infty  {\frac{{{{\left( {{\gamma _1}} \right)}_{{\alpha _1}m + {\beta _1}n}}{{\left( {{\gamma _2}} \right)}_{{\alpha _2}m}}{{\left( {{\gamma _3}} \right)}_{{\beta _2}n}}}}{{\Gamma \left( {{\delta _1} + {\alpha _3}m} \right)\Gamma \left( {{\delta _2} + {\beta _3}n} \right)}}\frac{{{x^m}}}{{\Gamma \left( {{\delta _3} + {\alpha _4}m} \right)}}\frac{{{y^n}}}{{\Gamma \left( {{\delta _4} + {\beta _4}n} \right)}}} ,
\end{array}
\end{equation}
\begin{equation}\label{eq 2.3}
\begin{array}{l}
{D_3}\left( {\begin{array}{*{20}{c}}
{{\gamma _1},{\alpha _1};{\gamma _2},{\beta _1};{\gamma _3},{\alpha _2};{\gamma _4},{\beta _2};}\\
{{\delta _1},{\alpha _3},{\beta _3};{\delta _2},{\alpha _4};{\delta _3},{\beta _4};}
\end{array}\left| {\begin{array}{*{20}{c}}
x\\
y
\end{array}} \right.} \right) \\
= \displaystyle \sum\limits_{m,n = 0}^\infty  {\frac{{{{\left( {{\gamma _1}} \right)}_{{\alpha _1}m}}{{\left( {{\gamma _2}} \right)}_{{\beta _1}n}}{{\left( {{\gamma _3}} \right)}_{{\alpha _2}m}}{{\left( {{\gamma _4}} \right)}_{{\beta _2}n}}}}{{\Gamma \left( {{\delta _1} + {\alpha _3}m + {\beta _3}n} \right)}}\frac{{{x^m}}}{{\Gamma \left( {{\delta _2} + {\alpha _4}m} \right)}}\frac{{{y^n}}}{{\Gamma \left( {{\delta _3} + {\beta _4}n} \right)}},}
\end{array}
\end{equation}
\begin{equation}\label{eq 2.4}
\begin{array}{l}
{D_4}\left( {\begin{array}{*{20}{c}}
{{\gamma _1},{\alpha _1},{\beta _1};{\gamma _2},{\alpha _2},{\beta _2};}\\
{{\delta _1},{\alpha _3};{\delta _2},{\beta _3};{\delta _3},{\alpha _4};{\delta _4},{\beta _4};}
\end{array}\left| {\begin{array}{*{20}{c}}
x\\
y
\end{array}} \right.} \right)\\
 =\displaystyle \sum\limits_{m,n = 0}^\infty  {\frac{{{{\left( {{\gamma _1}} \right)}_{{\alpha _1}m + {\beta _1}n}}{{\left( {{\gamma _2}} \right)}_{{\alpha _2}m + {\beta _2}n}}}}{{\Gamma \left( {{\delta _1} + {\alpha _3}m} \right)\Gamma \left( {{\delta _2} + {\beta _3}n} \right)}}\frac{{{x^m}}}{{\Gamma \left( {{\delta _3} + {\alpha _4}m} \right)}}\frac{{{y^n}}}{{\Gamma \left( {{\delta _4} + {\beta _4}n} \right)}}} ,
\end{array}
\end{equation}
\begin{equation}\label{eq 2.5}
\begin{array}{l}
{D_5}\left( {\begin{array}{*{20}{c}}
{{\gamma _1},{\alpha _1},{\beta _1};{\gamma _2},{\alpha _2},{\beta _2};}\\
{{\delta _1},{\alpha _3},{\beta _3};{\delta _2},{\alpha _4};{\delta _3},{\beta _4};}
\end{array}\left| {\begin{array}{*{20}{c}}
x\\
y
\end{array}} \right.} \right) \\
=\displaystyle \sum\limits_{m,n = 0}^\infty  {\frac{{{{\left( {{\gamma _1}} \right)}_{{\alpha _1}m + {\beta _1}n}}{{\left( {{\gamma _2}} \right)}_{{\alpha _2}m + {\beta _2}n}}}}{{\Gamma \left( {{\delta _1} + {\alpha _3}m + {\beta _3}n} \right)}}\frac{{{x^m}}}{{\Gamma \left( {{\delta _2} + {\alpha _4}m} \right)}}\frac{{{y^n}}}{{\Gamma \left( {{\delta _3} + {\beta _4}n} \right)}}} ,
\end{array}
\end{equation}

\begin{equation}\label{eq 2.6}
\begin{array}{l}
{E_3}\left( {\begin{array}{*{20}{c}}
{{\gamma _1},{\alpha _1},{\beta _1};{\gamma _2},{\alpha _2};}\\
{{\delta _1},{\alpha _3};{\delta _2},{\beta _2};{\delta _3},{\alpha _4};{\delta _4},{\beta _3};}
\end{array}\left| {\begin{array}{*{20}{c}}
x\\
y
\end{array}} \right.} \right)\\
 = \displaystyle \sum\limits_{m,n = 0}^\infty  {\frac{{{{\left( {{\gamma _1}} \right)}_{{\alpha _1}m + {\beta _1}n}}{{\left( {{\gamma _2}} \right)}_{{\alpha _2}m}}}}{{\Gamma \left( {{\delta _1} + {\alpha _3}m} \right)\Gamma \left( {{\delta _2} + {\beta _2}n} \right)}}\frac{{{x^m}}}{{\Gamma \left( {{\delta _3} + {\alpha _4}m} \right)}}\frac{{{y^n}}}{{\Gamma \left( {{\delta _4} + {\beta _3}n} \right)}}} ,
\end{array}
\end{equation}
\begin{equation}\label{eq 2.7}
\begin{array}{l}
{E_4}\left( {\begin{array}{*{20}{c}}
{{\gamma _1},{\alpha _1},{\beta _1};}\\
{{\delta _1},{\alpha _2};{\delta _2},{\beta _2};{\delta _3},{\alpha _3};{\delta _4},{\beta _3};}
\end{array}\left| {\begin{array}{*{20}{c}}
x\\
y
\end{array}} \right.} \right)\\
 = \displaystyle \sum\limits_{m,n = 0}^\infty  {\frac{{{{\left( {{\gamma _1}} \right)}_{{\alpha _1}m + {\beta _1}n}}}}{{\Gamma \left( {{\delta _1} + {\alpha _2}m} \right)\Gamma \left( {{\delta _2} + {\beta _2}n} \right)}}\frac{{{x^m}}}{{\Gamma \left( {{\delta _3} + {\alpha _3}m} \right)}}\frac{{{y^n}}}{{\Gamma \left( {{\delta _4} + {\beta _3}n} \right)}}} ,
\end{array}
\end{equation}
\begin{equation}\label{eq 2.8}
\begin{array}{l}
{E_5}\left( {\begin{array}{*{20}{c}}
{{\gamma _1},{\alpha _1};}\\
{{\delta _1},{\alpha _2},{\beta _1};{\delta _2},{\alpha _3};{\delta _3},{\beta _2};}
\end{array}\left| {\begin{array}{*{20}{c}}
x\\
y
\end{array}} \right.} \right) \\
= \displaystyle \sum\limits_{m,n = 0}^\infty  {\frac{{{{\left( {{\gamma _1}} \right)}_{{\alpha _1}m}}}}{{\Gamma \left( {{\delta _1} + {\alpha _2}m + {\beta _1}n} \right)}}\frac{{{x^m}}}{{\Gamma \left( {{\delta _2} + {\alpha _3}m} \right)}}\frac{{{y^n}}}{{\Gamma \left( {{\delta _3} + {\beta _2}n} \right)}},}
\end{array}
\end{equation}
\begin{equation}\label{eq 2.9}
\begin{array}{l}
{E_6}\left( {\begin{array}{*{20}{c}}
{{\gamma _1},{\alpha _1};{\gamma _2},{\beta _1};{\gamma _3},{\alpha _2};}\\
{{\delta _1},{\alpha _3},{\beta _2};{\delta _2},{\alpha _4};{\delta _3},{\beta _3};}
\end{array}\left| {\begin{array}{*{20}{c}}
x\\
y
\end{array}} \right.} \right) \\
= \displaystyle \sum\limits_{m,n = 0}^\infty  {\frac{{{{\left( {{\gamma _1}} \right)}_{{\alpha _1}m}}{{\left( {{\gamma _2}} \right)}_{{\beta _1}n}}{{\left( {{\gamma _3}} \right)}_{{\alpha _2}m}}}}{{\Gamma \left( {{\delta _1} + {\alpha _3}m + {\beta _2}n} \right)}}\frac{{{x^m}}}{{\Gamma \left( {{\delta _2} + {\alpha _4}m} \right)}}\frac{{{y^n}}}{{\Gamma \left( {{\delta _3} + {\beta _3}n} \right)}}} ,
\end{array}
\end{equation}
\begin{equation}\label{eq 2.10}
\begin{array}{l}
{E_7}\left( {\begin{array}{*{20}{c}}
{{\gamma _1},{\alpha _1};{\gamma _2},{\alpha _2};}\\
{{\delta _1},{\alpha _3},{\beta _1};{\delta _2},{\alpha _4};{\delta _3},{\beta _2};}
\end{array}\left| {\begin{array}{*{20}{c}}
x\\
y
\end{array}} \right.} \right) \\
= \displaystyle \sum\limits_{m,n = 0}^\infty  {\frac{{{{\left( {{\gamma _1}} \right)}_{{\alpha _1}m}}{{\left( {{\gamma _2}} \right)}_{{\alpha _2}m}}}}{{\Gamma \left( {{\delta _1} + {\alpha _3}m + {\beta _1}n} \right)}}\frac{{{x^m}}}{{\Gamma \left( {{\delta _2} + {\alpha _4}m} \right)}}\frac{{{y^n}}}{{\Gamma \left( {{\delta _3} + {\beta _2}n} \right)}}} ,
\end{array}
\end{equation}
\begin{equation}\label{eq 2.11}
\begin{array}{l}
{E_8}\left( {\begin{array}{*{20}{c}}
{{\gamma _1},{\alpha _1},{\beta _1};}\\
{{\delta _1},{\alpha _2},{\beta _2};{\delta _2},{\alpha _3};{\delta _3},{\beta _3};}
\end{array}\left| {\begin{array}{*{20}{c}}
x\\
y
\end{array}} \right.} \right)\\
 = \displaystyle \sum\limits_{m,n = 0}^\infty  {\frac{{{{\left( {{\gamma _1}} \right)}_{{\alpha _1}m + {\beta _1}n}}}}{{\Gamma \left( {{\delta _1} + {\alpha _2}m + {\beta _2}n} \right)}}\frac{{{x^m}}}{{\Gamma \left( {{\delta _2} + {\alpha _3}m} \right)}}\frac{{{y^n}}}{{\Gamma \left( {{\delta _3} + {\beta _3}n} \right)}}} ,
\end{array}
\end{equation}
\begin{equation}\label{eq 2.12}
\begin{array}{l}
{E_9}\left( {\begin{array}{*{20}{c}}
{ - ;}\\
{{\delta _1},{\alpha _1},{\beta _1};{\delta _2},{\alpha _2};{\delta _3},{\beta _2};}
\end{array}\left| {\begin{array}{*{20}{c}}
x\\
y
\end{array}} \right.} \right) \\
=  \displaystyle  \sum\limits_{m,n = 0}^\infty  {\frac{1}{{\Gamma \left( {{\delta _1} + {\alpha _1}m + {\beta _1}n} \right)}}\frac{{{x^m}}}{{\Gamma \left( {{\delta _2} + {\alpha _2}m} \right)}}\frac{{{y^n}}}{{\Gamma \left( {{\delta _3} + {\beta _2}n} \right)}},}
\end{array}
\end{equation}
\begin{equation}\label{eq 2.13}
\begin{array}{l}
{E_{10}}\left( {\begin{array}{*{20}{c}}
{{\gamma _1},{\alpha _1},{\beta _1};{\gamma _2},{\alpha _2};}\\
{{\delta _1},{\alpha _3};{\delta _2},{\beta _2};}
\end{array}\left| {\begin{array}{*{20}{c}}
x\\
y
\end{array}} \right.} \right) \\
=\displaystyle \sum\limits_{m,n = 0}^\infty  {{{\left( {{\gamma _1}} \right)}_{{\alpha _1}m + {\beta _1}n}}{{\left( {{\gamma _2}} \right)}_{{\alpha _2}m}}\frac{{{x^m}}}{{\Gamma \left( {{\delta _1} + {\alpha _3}m} \right)}}\frac{{{y^n}}}{{\Gamma \left( {{\delta _2} + {\beta _2}n} \right)}}.}
\end{array}
\end{equation}
\begin{equation}\label{eq 2.14}
\begin{array}{l}
    {{E}_{11}}\left( \begin{matrix}
   {{\gamma }_{1}},{{\alpha }_{1}},{{\beta }_{1}};{{\gamma }_{2}},{{\alpha }_{2}};{{\gamma }_{3}},{{\beta }_{2}};  \\
   {{\delta }_{1}},{{\alpha }_{3}};{{\delta }_{2}},{{\beta }_{3}};  \\
\end{matrix}\left| \begin{matrix}
   x  \\
   y  \\
\end{matrix} \right. \right) \\
 = \displaystyle \sum\limits_{m,n=0}^{\infty }{{{\left( {{\gamma }_{1}} \right)}_{{{\alpha }_{1}}m+{{\beta }_{1}}n}}{{\left( {{\gamma }_{2}} \right)}_{{{\alpha }_{2}}m}}{{\left( {{\gamma }_{3}} \right)}_{{{\beta }_{2}}n}}\frac{{{x}^{m}}}{\Gamma \left( {{\delta }_{1}}+{{\alpha }_{3}}m \right)}\frac{{{y}^{n}}}{\Gamma \left( {{\delta }_{2}}+{{\beta }_{3}}n \right)}}. \\
\end{array}
\end{equation}
Note that the introduced generalized Mittag-Leffler functions (\ref{eq 2.1}) - (\ref{eq 2.14}) in particular values of the parameters coincide with the known hypergeometric functions. For example, if in (\ref{eq 2.1}) the parameters take the following values ${\alpha _1} = {\alpha _2} = {\alpha _3} = {\alpha _4} = {\beta _1} = {\beta _2} = {\beta _3} = {\beta _4} = {\delta _2} = {\delta _3} = 1,$  it coincides with the Appel function, that is ${D_1}(x,y) = {F_1}(x,y)$  \cite{42}.

\section{Region of the convergence}

Following Horn \cite{2}, we determine the region of convergence of the introduced Mittag-Leffler-type function $D_1$.

{\bf Definition.}
Let us call positive values $r$, $s$   the associated radii of convergence of the double series
\begin{equation}\label{eq 18}
\sum\limits_{m,n = 0}^\infty  {A(m,n){x^m}{y^n}},
\end{equation}
if it converges absolutely at $\left| x \right| < r $,   $\left| y \right| < s $  and diverges at $\left| x \right| > r $,   $\left| y \right| > s $.

Let us also assume that the $\max (r) = R,$  $\max (s) = S.$  Points $(r,s)$  lie on the curve $C$, which is located entirely in the rectangle $0 < r < R,\,0 < s < S.$ This curve divides the rectangle into two parts; the part containing the point $r=s=0$  is a two-dimensional image of the region of convergence of the double power series. Studying the convergence of the series (\ref{eq 18}), Horn introduced the functions
\begin{equation*}
\Phi \left( {\mu ,\nu } \right) = \mathop {\lim }\limits_{t \to \infty } f\left( {\mu t,\nu t} \right),\,\,\,\Psi \left( {\mu ,\nu } \right) = \mathop {\lim }\limits_{t \to \infty } g\left( {\mu t,\nu t} \right),
\end{equation*}
where
\begin{equation}\label{eq 19}
f\left( {m,n} \right) = \frac{{A(m + 1,n)}}{{A(m,n)}},\,\,\,g\left( {m,n} \right) = \frac{{A(m,n + 1)}}{{A(m,n)}},
\end{equation}
and showed that $R = {\left| {\Phi \left( {1,0} \right)} \right|^{ - 1}},\,S = {\left| {\Psi \left( {0,1} \right)} \right|^{ - 1}}$ and that $C$ has a parametric representation $r = {\left| {\Phi \left( {\mu ,\nu } \right)} \right|^{ - 1}},\,\,s = {\left| {\Psi \left( {\mu ,\nu } \right)} \right|^{ - 1}},\,\,\mu ,\,\,\nu  > 0.$\\
Consider the function (\ref{eq 19}). It follows from the definition of the function $D_1$
\begin{equation*}
\begin{array}{l}
f(\mu t,\nu t) = \displaystyle \frac{{\Gamma \left( {{\gamma _1} + {\alpha _1} + {\alpha _1}\mu t + {\beta _1}\nu t} \right)\Gamma \left( {{\gamma _2} + {\alpha _2} + {\alpha _2}\mu t} \right)}}{{\Gamma \left( {{\gamma _1} + {\alpha _1}\mu t + {\beta _1}\nu t} \right)\Gamma \left( {{\gamma _2} + {\alpha _2}\mu t} \right)}}\\
\,\,\,\,\,\,\,\,\,\,\,\,\,\,\,\,\,\,\,\,\, \times \displaystyle \frac{{\Gamma \left( {{\delta _1} + {\alpha _3}\mu t + {\beta _3}\nu t} \right)\Gamma \left( {{\delta _2} + {\alpha _4}\mu t} \right)}}{{\Gamma \left( {{\delta _1} + {\alpha _3} + {\alpha _3}\mu t + {\beta _3}\nu t} \right)\Gamma \left( {{\delta _2} + {\alpha _4} + {\alpha _4}\mu t} \right)}},\\
g(\mu t,\nu t) = \displaystyle \frac{{\Gamma \left( {{\gamma _1} + {\beta _1} + {\alpha _1}\mu t + {\beta _1}\nu t} \right)\Gamma \left( {{\gamma _3} + {\beta _2} + {\beta _2}\nu t} \right)}}{{\Gamma \left( {{\gamma _1} + {\alpha _1}\mu t + {\beta _1}\nu t} \right)\Gamma \left( {{\gamma _3} + {\beta _2}\nu t} \right)}}\\
\,\,\,\,\,\,\,\,\,\,\,\,\,\,\,\,\,\,\,\, \times \displaystyle \frac{{\Gamma \left( {{\delta _1} + {\alpha _3}\mu t + {\beta _3}\nu t} \right)\Gamma \left( {{\delta _3} + {\beta _4}\nu t} \right)}}{{\Gamma \left( {{\delta _1} + {\beta _3} + {\alpha _3}\mu t + {\beta _3}\nu t} \right)\Gamma \left( {{\delta _3} + {\beta _4} + {\beta _4}\nu t} \right)}}.
\end{array}
\end{equation*}
Due to the asymptotics of the Gamma function for large arguments \cite{43}
\begin{equation*}
\frac{{\Gamma \left( {z + \alpha } \right)}}{{\Gamma \left( {z + \beta } \right)}} \sim {z^{\alpha  - \beta }}\left[ {1 + \frac{{\left( {\alpha  - \beta } \right)\left( {\alpha  + \beta  - 1} \right)}}{{2z}} + O\left( {{z^{ - 2}}} \right)} \right],\,\,\left| {\arg \left( z \right)} \right| \le \pi,
\end{equation*}
we have
\begin{equation}\label{eq 20}
\begin{array}{l}
f(\mu t,\nu t) \sim \displaystyle\frac{1}{E}{t^{ - \Delta }},\,\,\,\Delta  = {\alpha _3} + {\alpha _4} - {\alpha _1} - {\alpha _2},\\
E = \displaystyle {\mu ^{{\alpha _4} - {\alpha _2}}}\frac{{{{\left( {{\alpha _3}\mu  + {\beta _3}\nu } \right)}^{{\alpha _3}}}{{\left( {{\alpha _4}} \right)}^{{\alpha _4}}}}}{{{{\left( {{\alpha _1}\mu  + {\beta _1}\nu } \right)}^{{\alpha _1}}}{{\left( {{\alpha _2}} \right)}^{{\alpha _2}}}}},\,\,\,\,\,G = \frac{{{{\left( {{\alpha _3}} \right)}^{{\alpha _3}}}{{\left( {{\alpha _4}} \right)}^{{\alpha _4}}}}}{{{{\left( {{\alpha _1}} \right)}^{{\alpha _1}}}{{\left( {{\alpha _2}} \right)}^{{\alpha _2}}}}},
\end{array}
\end{equation}
Similarly, we define
\begin{equation*}
\begin{array}{l}
g(\mu t,\nu t) \sim \displaystyle\frac{1}{{E'}} \cdot {t^{ - \Delta '}},\,\,\Delta ' = {\beta _3} + {\beta _4} - {\beta _1} - {\beta _2},\\
E' = \displaystyle {\nu ^{{\beta _4} - {\beta _2}}}\frac{{{{\left( {{\alpha _3}\mu  + {\beta _3}\nu } \right)}^{{\beta _3}}}{{\left( {{\beta _4}} \right)}^{{\beta _4}}}}}{{{{\left( {{\alpha _1}\mu  + {\beta _1}\nu } \right)}^{{\beta _1}}}{{\left( {{\beta _2}} \right)}^{{\beta _2}}}}},\,\,G' = \frac{{{{\left( {{\beta _3}} \right)}^{{\beta _3}}}{{\left( {{\beta _4}} \right)}^{{\beta _4}}}}}{{{{\left( {{\beta _1}} \right)}^{{\beta _1}}}{{\left( {{\beta _2}} \right)}^{{\beta _2}}}}}.
\end{array}
\end{equation*}
Now consider some cases:\\
{\bf Case 1.} Let $\Delta  > 0,\,\,\,\Delta ' > 0.$ Then from (\ref{eq 19}) and (\ref{eq 20}) it follows
\begin{equation*}
\begin{array}{l}
\Phi \left( {\mu ,\nu } \right) = \mathop {\lim }\limits_{t \to \infty } f(\mu t,\nu t) = 0,\,\,\Psi \left( {\mu ,\nu } \right) = \mathop {\lim }\limits_{t \to \infty } g(\mu t,\nu t) = 0.
\end{array}
\end{equation*}
Positive numbers $r$  and $s$   are large numbers. The series converges for any value of the argument.\\
{\bf Case 2.} Let $\Delta  = 0,\,\,\,\Delta ' = 0.$ Then from (\ref{eq 19}) and (\ref{eq 20}) it follows
\begin{equation*}
\begin{array}{l}
\Phi \left( {\mu ,\nu } \right) = \displaystyle \frac{1}{E},\,\,\,\,\,\,\,\,\Psi \left( {\mu ,\nu } \right) = \displaystyle \frac{1}{{E'}},
\end{array}
\end{equation*}
which immediately led us to the parametric representation of the curve $C$ on the plane  $(r,s)$ in the form $r = G,\,\,\,\,s = G'.$ Therefore, the series converges absolutely for the values of  $\left| x \right| < \rho $  and $\left| y \right| < \rho '$ where $\rho  = \mathop {\min }\limits_{\mu ,\nu  > 0} (E),\,\,\,\rho ' = \mathop {\min }\limits_{\mu ,\nu  > 0} (E').$ \\
{\bf Case 3.} Let $\Delta  < 0,\,\,\,\Delta ' < 0.$ The series diverges $r = s = 0.$  The series converges only at the point $x = y = 0.$ \\
{\bf Case 4.} Let $\Delta  = 0,\,\,\,\Delta ' > 0.$ Then the series converges absolutely in the region $\left| x \right| < \rho$  and  $\left| y \right| < \infty$.\\
{\bf Case 5.} Let $\Delta  > 0,\,\,\,\Delta ' = 0.$ Then the series converges absolutely in the region $\left| x \right| < \infty$  and $\left| y \right| < \rho'$.

\section{Integral representations}

For a generalized hypergeometric function of the Mittag-Leffler-type $D_1$, the following integral representations of the Euler type are valid:
\begin{equation*}
\begin{array}{l}
  {{D}_{1}}\left( \begin{matrix}
   {{\gamma }_{1}},{{\alpha }_{1}},{{\beta }_{1}};{{\gamma }_{2}},{{\alpha }_{2}};{{\gamma }_{3}},{{\beta }_{2}};  \\
   {{\delta }_{1}},{{\alpha }_{3}},{{\beta }_{3}};{{\delta }_{2}},{{\alpha }_{4}};{{\delta }_{3}},{{\beta }_{4}};  \\
\end{matrix}\left| \begin{matrix}
   x  \\
   y  \\
\end{matrix} \right. \right) =\displaystyle \frac{\Gamma \left( \mu  \right)}{\Gamma \left( {{\gamma }_{2}} \right)\Gamma \left( \mu -{{\gamma }_{2}} \right)}\\
  \times \displaystyle \int\limits_{0}^{1}{{}}{{\xi }^{{{\gamma }_{2}}-1}}{{\left( 1-\xi  \right)}^{\mu -{{\gamma }_{2}}-1}}{{D}_{1}}\left( \begin{matrix}
   {{\gamma }_{1}},{{\alpha }_{1}},{{\beta }_{1}};\mu ,{{\alpha }_{2}};{{\gamma }_{3}},{{\beta }_{2}};  \\
   {{\delta }_{1}},{{\alpha }_{3}},{{\beta }_{3}};{{\delta }_{2}},{{\alpha }_{4}};{{\delta }_{3}},{{\beta }_{4}};  \\
\end{matrix}\left| \begin{matrix}
   x{{\xi }^{{{\alpha }_{2}}}}  \\
   y  \\
\end{matrix} \right. \right)d\xi , \\
 \operatorname{Re}\mu >\operatorname{Re}{{\gamma }_{2}}>0, \\
\end{array}
\end{equation*}
\begin{equation*}
\begin{array}{l}
  {{D}_{1}}\left( \begin{matrix}
   {{\gamma }_{1}},{{\alpha }_{1}},{{\beta }_{1}};{{\gamma }_{2}},{{\alpha }_{2}};{{\gamma }_{3}},{{\beta }_{2}};  \\
   {{\delta }_{1}},{{\alpha }_{3}},{{\beta }_{3}};{{\delta }_{2}},{{\alpha }_{4}};{{\delta }_{3}},{{\beta }_{4}};  \\
\end{matrix}\left| \begin{matrix}
   x  \\
   y  \\
\end{matrix} \right. \right)= \displaystyle \frac{\Gamma \left( \mu  \right)}{\Gamma \left( {{\gamma }_{3}} \right)\Gamma \left( \mu -{{\gamma }_{3}} \right)} \\
 \times \int\limits_{0}^{1}{{}}{{\xi }^{{{\gamma }_{3}}-1}}{{\left( 1-\xi  \right)}^{\mu -{{\gamma }_{3}}-1}}{{D}_{1}}\left( \begin{matrix}
   {{\gamma }_{1}},{{\alpha }_{1}},{{\beta }_{1}};{{\gamma }_{2}},{{\alpha }_{2}};\mu ,{{\beta }_{2}};  \\
   {{\delta }_{1}},{{\alpha }_{3}},{{\beta }_{3}};{{\delta }_{2}},{{\alpha }_{4}};{{\delta }_{3}},{{\beta }_{4}};  \\
\end{matrix}\left| \begin{matrix}
   x  \\
   y{{\xi }^{{{\beta }_{2}}}}  \\
\end{matrix} \right. \right)d\xi , \\
 \operatorname{Re}\mu >\operatorname{Re}{{\gamma }_{3}}>0, \\
\end{array}
\end{equation*}
\begin{equation*}
\begin{array}{l}
  {{D}_{1}}\left( \begin{matrix}
   {{\gamma }_{1}},{{\alpha }_{1}},{{\beta }_{1}};{{\gamma }_{2}},{{\alpha }_{2}};{{\gamma }_{3}},{{\beta }_{2}};  \\
   {{\delta }_{1}},{{\alpha }_{3}},{{\beta }_{3}};{{\delta }_{2}},{{\alpha }_{4}};{{\delta }_{3}},{{\beta }_{4}};  \\
\end{matrix}\left| \begin{matrix}
   x  \\
   y  \\
\end{matrix} \right. \right)= \displaystyle \frac{\Gamma \left( {{\mu }_{1}} \right)\Gamma \left( {{\mu }_{2}} \right)}{\Gamma \left( {{\gamma }_{2}} \right)\Gamma \left( {{\gamma }_{3}} \right)\Gamma \left( {{\mu }_{1}}-{{\gamma }_{2}} \right)\Gamma \left( {{\mu }_{2}}-{{\gamma }_{3}} \right)} \\
 \times \displaystyle \int\limits_{0}^{1}{\int\limits_{0}^{1}{{}}}{{\xi }^{{{\gamma }_{2}}-1}}{{\eta }^{{{\gamma }_{3}}-1}}{{\left( 1-\xi  \right)}^{{{\mu }_{1}}-{{\gamma }_{2}}-1}}{{\left( 1-\eta  \right)}^{{{\mu }_{2}}-{{\gamma }_{3}}-1}} \\ \times {{D}_{1}}\left( \begin{matrix}
   {{\gamma }_{1}},{{\alpha }_{1}},{{\beta }_{1}};{{\mu }_{1}},{{\alpha }_{2}};{{\mu }_{2}},{{\beta }_{2}};  \\
   {{\delta }_{1}},{{\alpha }_{3}},{{\beta }_{3}};{{\delta }_{2}},{{\alpha }_{4}};{{\delta }_{3}},{{\beta }_{4}};  \\
\end{matrix}\left| \begin{matrix}
   x{{\xi }^{{{\alpha }_{2}}}}  \\
   y{{\eta }^{{{\beta }_{2}}}}  \\
\end{matrix} \right. \right)d\xi d\eta,\\
 \operatorname{Re}{{\mu }_{1}}>\operatorname{Re}{{\gamma }_{2}}>0,\,\,\operatorname{Re}{{\mu }_{2}}>\operatorname{Re}{{\gamma }_{3}}>0, \\
\end{array}
\end{equation*}
\begin{equation*}
\begin{array}{l}
   \displaystyle \int\limits_{0}^{1}{\int\limits_{0}^{1}{{}}}{{\xi }^{{{\delta }_{2}}-1}}{{\eta }^{{{\delta }_{3}}-1}}{{\left( 1-\xi  \right)}^{{{\sigma }_{1}}-1}}{{\left( 1-\eta  \right)}^{{{\sigma }_{2}}-1}}{{D}_{1}}\left( \begin{matrix}
   {{\gamma }_{1}},{{\alpha }_{1}},{{\beta }_{1}};{{\gamma }_{2}},{{\alpha }_{2}};{{\gamma }_{3}},{{\beta }_{2}};  \\
   {{\delta }_{1}},{{\alpha }_{3}},{{\beta }_{3}};{{\delta }_{2}},{{\alpha }_{4}};{{\delta }_{3}},{{\beta }_{4}};  \\
\end{matrix}\left| \begin{matrix}
   x{{\xi }^{{{\alpha }_{4}}}}  \\
   y{{\eta }^{{{\beta }_{4}}}}  \\
\end{matrix} \right. \right)d\xi d\eta = \\
 =\Gamma \left( {{\sigma }_{1}} \right)\Gamma \left( {{\sigma }_{2}} \right){{D}_{1}}\left( \begin{matrix}
   {{\gamma }_{1}},{{\alpha }_{1}},{{\beta }_{1}};{{\gamma }_{2}},{{\alpha }_{2}};{{\gamma }_{3}},{{\beta }_{2}};  \\
   {{\delta }_{1}},{{\alpha }_{3}},{{\beta }_{3}};{{\delta }_{2}}+{{\sigma }_{1}},{{\alpha }_{4}};{{\delta }_{3}}+{{\sigma }_{2}},{{\beta }_{4}};  \\
\end{matrix}\left| \begin{matrix}
   x  \\
   y  \\
\end{matrix} \right. \right),\\
\operatorname{Re}{{\sigma }_{1}}>0,\,\,\operatorname{Re}{{\sigma }_{2}}>0, \\
\end{array}
\end{equation*}
\begin{equation*}
\begin{array}{l}
   {{D}_{1}}\left( \begin{matrix}
   {{\gamma }_{1}},{{\alpha }_{1}},{{\beta }_{1}};{{\gamma }_{2}},{{\alpha }_{2}};{{\gamma }_{3}},{{\beta }_{2}};  \\
   {{\delta }_{1}},{{\alpha }_{3}},{{\beta }_{3}};{{\delta }_{2}},{{\alpha }_{4}};{{\delta }_{3}},{{\beta }_{4}};  \\
\end{matrix}\left| \begin{matrix}
   x  \\
   y  \\
\end{matrix} \right. \right)= \displaystyle \frac{\Gamma \left( \mu  \right)}{\Gamma \left( {{\gamma }_{1}} \right)\Gamma \left( \mu -{{\gamma }_{1}} \right)}  \\
 \times \displaystyle  \int\limits_{0}^{1}{{}}{{\xi }^{{{\gamma }_{1}}-1}}{{\left( 1-\xi  \right)}^{\mu -{{\gamma }_{1}}-1}}{{D}_{1}}\left( \begin{matrix}
   \mu ,{{\alpha }_{1}},{{\beta }_{1}};{{\gamma }_{2}},{{\alpha }_{2}};{{\gamma }_{3}},{{\beta }_{2}};  \\
   {{\delta }_{1}},{{\alpha }_{3}},{{\beta }_{3}};{{\delta }_{2}},{{\alpha }_{4}};{{\delta }_{3}},{{\beta }_{4}};  \\
\end{matrix}\left| \begin{matrix}
   x{{\xi }^{{{\alpha }_{1}}}}  \\
   y{{\xi }^{{{\beta }_{1}}}}  \\
\end{matrix} \right. \right)d\xi ,\,\,\operatorname{Re}\mu >\operatorname{Re}{{\gamma }_{1}}>0, \\
\end{array}
\end{equation*}
\begin{equation*}
\begin{array}{l}
 {{D}_{1}}\left( \begin{matrix}
   {{\gamma }_{1}},{{\alpha }_{1}},{{\beta }_{1}};{{\gamma }_{2}},{{\alpha }_{2}};{{\gamma }_{3}},{{\beta }_{2}};  \\
   {{\mu }_{1}}+{{\mu }_{2}},{{\alpha }_{3}},{{\beta }_{3}};{{\delta }_{2}},{{\alpha }_{4}};{{\delta }_{3}},{{\beta }_{4}};  \\
\end{matrix}\left| \begin{matrix}
   x  \\
   y  \\
\end{matrix} \right. \right) \\
 = \displaystyle \int\limits_{0}^{1}{{}}{{\xi }^{{{\mu }_{1}}-1}}{{\left( 1-\xi  \right)}^{{{\mu }_{2}}-1}}{{D}_{2}}\left( \begin{matrix}
   {{\gamma }_{1}},{{\alpha }_{1}},{{\beta }_{1}};{{\gamma }_{2}},{{\alpha }_{2}};{{\gamma }_{3}},{{\beta }_{2}};  \\
   {{\mu }_{1}},{{\alpha }_{3}};{{\mu }_{2}},{{\beta }_{3}};{{\delta }_{2}},{{\alpha }_{4}};{{\delta }_{3}},{{\beta }_{4}};  \\
\end{matrix}\left| \begin{matrix}
   x{{\xi }^{{{\alpha }_{3}}}}  \\
   y{{\left( 1-\xi  \right)}^{{{\beta }_{3}}}}  \\
\end{matrix} \right. \right)d\xi , \\
 \operatorname{Re}{{\mu }_{1}}>0,\operatorname{Re}{{\mu }_{2}}>0, \\
\end{array}
\end{equation*}
\begin{equation*}
\begin{array}{l}
  {{D}_{1}}\left( \begin{matrix}
   {{\gamma }_{1}},{{\alpha }_{1}},{{\beta }_{1}};{{\gamma }_{2}},{{\alpha }_{2}};{{\gamma }_{3}},{{\beta }_{2}};  \\
   {{\delta }_{1}},{{\alpha }_{1}},{{\beta }_{1}};{{\delta }_{2}},{{\alpha }_{4}};{{\delta }_{3}},{{\beta }_{4}};  \\
\end{matrix}\left| \begin{matrix}
   x  \\
   y  \\
\end{matrix} \right. \right)= \displaystyle \frac{1}{\Gamma \left( {{\gamma }_{1}} \right)\Gamma \left( {{\delta }_{1}}-{{\gamma }_{1}} \right)}\\
\times \displaystyle \int\limits_{0}^{1}{{}}{{\xi }^{{{\gamma }_{1}}-1}}{{\left( 1-\xi  \right)}^{{{\delta }_{1}}-{{\gamma }_{1}}-1}}E\left( \begin{matrix}
   {{\gamma }_{2}},{{\alpha }_{2}};  \\
   {{\delta }_{2}},{{\alpha }_{4}};  \\
\end{matrix}x{{\xi }^{{{\alpha }_{1}}}} \right)E\left( \begin{matrix}
   {{\gamma }_{3}},{{\beta }_{2}};  \\
   {{\delta }_{3}},{{\beta }_{4}};  \\
\end{matrix}y{{\xi }^{{{\beta }_{1}}}} \right)d\xi , \\
  \operatorname{Re}{{\delta }_{1}}>\operatorname{Re}{{\gamma }_{1}}>0, \\
\end{array}
\end{equation*}
\begin{equation*}
\begin{array}{l}
  {{D}_{1}}\left( \begin{matrix}
   {{\gamma }_{1}},{{\alpha }_{1}},{{\beta }_{1}};{{\gamma }_{2}},{{\alpha }_{2}};{{\gamma }_{3}},{{\beta }_{2}};  \\
   {{\delta }_{1}},{{\alpha }_{3}},{{\beta }_{3}};{{\delta }_{2}},{{\alpha }_{4}};{{\delta }_{3}},{{\beta }_{4}};  \\
\end{matrix}\left| \begin{matrix}
   x  \\
   y  \\
\end{matrix} \right. \right)= \displaystyle  \frac{\Gamma \left( {{\gamma }_{2}}+{{\gamma }_{3}} \right)}{\Gamma \left( {{\gamma }_{2}} \right)\Gamma \left( {{\gamma }_{3}} \right)} \\
 \times \displaystyle \int\limits_{0}^{1}{{}}{{\xi }^{{{\gamma }_{2}}-1}}{{\left( 1-\xi  \right)}^{{{\gamma }_{3}}-1}}{{D}_{5}}\left( \begin{matrix}
   {{\gamma }_{1}},{{\alpha }_{1}},{{\beta }_{1}};{{\gamma }_{2}}+{{\gamma }_{3}},{{\alpha }_{2}},{{\beta }_{2}};  \\
   {{\delta }_{1}},{{\alpha }_{3}},{{\beta }_{3}};{{\delta }_{2}},{{\alpha }_{4}};{{\delta }_{3}},{{\beta }_{4}};  \\
\end{matrix}\left| \begin{matrix}
   x{{\xi }^{{{\alpha }_{2}}}}  \\
   y{{\left( 1-\xi  \right)}^{{{\beta }_{2}}}}  \\
\end{matrix} \right. \right)d\xi \\
\operatorname{Re}{{\gamma }_{2}}>0,\,\,\operatorname{Re}{{\gamma }_{3}}>0, \\
\end{array}
\end{equation*}
\begin{equation*}
\begin{array}{l}
   {{D}_{1}}\left( \begin{matrix}
   {{\gamma }_{1}},{{\alpha }_{1}},{{\beta }_{1}};{{\gamma }_{2}},{{\alpha }_{2}};{{\gamma }_{3}},{{\beta }_{2}};  \\
   {{\delta }_{1}},{{\alpha }_{3}},{{\beta }_{3}};{{\delta }_{2}},{{\alpha }_{4}};{{\delta }_{3}},{{\beta }_{4}};  \\
\end{matrix}\left| \begin{matrix}
   x  \\
   y  \\
\end{matrix} \right. \right)= \displaystyle \frac{1}{\Gamma \left( {{\gamma }_{3}} \right)}\int\limits_{0}^{1}{{}}{{\xi }^{{{\gamma }_{3}}-1}}{{\left( 1-\xi  \right)}^{{{\delta }_{1}}-{{\gamma }_{3}}-1}} \\
  \times {{E}_{2}}\left( \begin{matrix}
   {{\gamma }_{1}},{{\alpha }_{1}},{{\beta }_{1}};{{\gamma }_{2}},{{\alpha }_{2}};  \\
   {{\delta }_{1}}-{{\gamma }_{3}},{{\alpha }_{3}},{{\beta }_{3}}-{{\beta }_{2}};{{\delta }_{2}},{{\alpha }_{4}};{{\delta }_{3}},{{\beta }_{3}};  \\
\end{matrix}\left| \begin{matrix}
   x{{\left( 1-\xi  \right)}^{{{\alpha }_{3}}}}  \\
   y{{\xi }^{{{\beta }_{2}}}}{{\left( 1-\xi  \right)}^{{{\beta }_{3}}-{{\beta }_{2}}}}  \\
\end{matrix} \right. \right)d\xi ,\,\,\\
\operatorname{Re}{{\delta }_{1}}>\operatorname{Re}{{\gamma }_{3}}>0, \\
\end{array}
\end{equation*}
\begin{equation}\label{eq 21}
\begin{array}{l}
    {{D}_{1}}\left( \begin{matrix}
   {{\gamma }_{1}},{{\alpha }_{3}}-{{\alpha }_{2}},{{\beta }_{1}};{{\gamma }_{2}},{{\alpha }_{2}};{{\gamma }_{3}},{{\beta }_{2}};  \\
   {{\delta }_{1}},{{\alpha }_{3}},{{\beta }_{1}};{{\delta }_{2}},{{\alpha }_{4}};{{\delta }_{3}},{{\beta }_{4}};  \\
\end{matrix}\left| \begin{matrix}
   x  \\
   y  \\
\end{matrix} \right. \right)= \displaystyle \frac{1}{\Gamma \left( {{\gamma }_{1}} \right)\Gamma \left( {{\gamma }_{2}} \right)\Gamma \left( {{\delta }_{1}}-{{\gamma }_{1}}-{{\gamma }_{2}} \right)}  \\
\times \displaystyle \int\limits_{0}^{1}{\int\limits_{0}^{1}{{}}}{{\xi }^{{{\gamma }_{2}}-1}}{{\eta }^{{{\gamma }_{1}}-1}}{{\left( 1-\xi  \right)}^{{{\delta }_{1}}-{{\gamma }_{2}}-1}}{{\left( 1-\eta  \right)}^{{{\delta }_{1}}-{{\gamma }_{1}}-{{\gamma }_{2}}-1}} \\
 \times E_{{{\delta }_{2}},{{\alpha }_{4}}}^{{{\gamma }_{3}},{{\beta }_{2}}}\left( x{{\xi }^{{{\alpha }_{2}}}}{{\left[ \left( 1-\xi  \right)\eta  \right]}^{\left( {{\alpha }_{3}}-{{\alpha }_{2}} \right)}} \right){{E}_{{{\delta }_{3}},{{\beta }_{4}}}}\left( y{{\left( 1-\xi  \right)}^{{{\beta }_{1}}}}{{\eta }^{{{\beta }_{1}}}} \right)d\xi d\eta ,\,\, \\
 \operatorname{Re}{{\gamma }_{1}}>0,\,\,\operatorname{Re}{{\gamma }_{2}}>0,\,\,\operatorname{Re}\left( {{\delta }_{1}}-{{\gamma }_{1}}-{{\gamma }_{2}} \right)>0,\,\,\,{{\alpha }_{3}}-{{\alpha }_{2}}>0. \\
\end{array}
\end{equation}
Let us prove the validity of the integral representation (\ref{eq 21}). 
\begin{proof}
We expand the integrands on the right-hand side into a series, then we have
\begin{equation*}
\begin{array}{l}
  I\left( x,y \right)= \displaystyle \frac{1}{\Gamma \left( {{\gamma }_{1}} \right)\Gamma \left( {{\gamma }_{2}} \right)\Gamma \left( {{\delta }_{1}}-{{\gamma }_{1}}-{{\gamma }_{2}} \right)}\times\\
  \int\limits_{0}^{1}{\int\limits_{0}^{1}{{}}}{{\xi }^{{{\gamma }_{2}}-1}}{{\eta }^{{{\gamma }_{1}}-1}}{{\left( 1-\xi  \right)}^{{{\delta }_{1}}-{{\gamma }_{2}}-1}}{{\left( 1-\eta  \right)}^{{{\delta }_{1}}-{{\gamma }_{1}}-{{\gamma }_{2}}-1}} \times\\
   \displaystyle \sum\limits_{m=0}^{\infty }{\frac{{{\left( {{\gamma }_{3}} \right)}_{{{\beta }_{2}}n}}{{\left( x{{\xi }^{{{\alpha }_{2}}}}{{\left[ \left( 1-\xi  \right)\eta  \right]}^{\left( {{\alpha }_{3}}-{{\alpha }_{2}} \right)}} \right)}^{m}}}{\Gamma \left( {{\delta }_{2}}+{{\alpha }_{4}}m \right)}}\sum\limits_{n=0}^{\infty }{\frac{{{\left( y{{\left( 1-\xi  \right)}^{{{\beta }_{1}}}}{{\eta }^{{{\beta }_{1}}}} \right)}^{n}}}{\Gamma \left( {{\delta }_{3}}+{{\beta }_{4}}n \right)}}d\xi d\eta  \\
  = \displaystyle \frac{1}{\Gamma \left( {{\gamma }_{1}} \right)\Gamma \left( {{\gamma }_{2}} \right)\Gamma \left( {{\delta }_{1}}-{{\gamma }_{1}}-{{\gamma }_{2}} \right)}\sum\limits_{m,n=0}^{\infty }{{}}\frac{{{\left( {{\gamma }_{3}} \right)}_{{{\beta }_{2}}n}}}{\Gamma \left( {{\delta }_{2}}+{{\alpha }_{4}}m \right)\Gamma \left( {{\delta }_{3}}+{{\beta }_{4}}n \right)}{{x}^{m}}{{y}^{n}} \\
  \times \displaystyle \int\limits_{0}^{1}{{}}{{\xi }^{{{\gamma }_{2}}+{{\alpha }_{2}}m-1}}{{\left( 1-\xi  \right)}^{{{\delta }_{1}}-{{\gamma }_{2}}+\left( {{\alpha }_{3}}-{{\alpha }_{2}} \right)m+{{\beta }_{1}}n-1}}d\xi \\
  \times \displaystyle \int\limits_{0}^{1}{{}}{{\eta }^{{{\gamma }_{1}}+\left( {{\alpha }_{3}}-{{\alpha }_{2}} \right)m+{{\beta }_{1}}n-1}}{{\left( 1-\eta  \right)}^{{{\delta }_{1}}-{{\gamma }_{1}}-{{\gamma }_{2}}-1}}d\eta = \\
  = \displaystyle \frac{1}{\Gamma \left( {{\gamma }_{1}} \right)\Gamma \left( {{\gamma }_{2}} \right)\Gamma \left( {{\delta }_{1}}-{{\gamma }_{1}}-{{\gamma }_{2}} \right)}\sum\limits_{m,n=0}^{\infty }{{}}\frac{{{\left( {{\gamma }_{3}} \right)}_{{{\beta }_{2}}n}}}{\Gamma \left( {{\delta }_{2}}+{{\alpha }_{4}}m \right)\Gamma \left( {{\delta }_{3}}+{{\beta }_{4}}n \right)}{{x}^{m}}{{y}^{n}} \\
  \times \displaystyle B\left( {{\gamma }_{2}}+{{\alpha }_{2}}m,{{\delta }_{1}}-{{\gamma }_{2}}+\left( {{\alpha }_{3}}-{{\alpha }_{2}} \right)m+{{\beta }_{1}}n \right)\times\\
  B\left( {{\gamma }_{1}}+\left( {{\alpha }_{3}}-{{\alpha }_{2}} \right)m+{{\beta }_{1}}n,{{\delta }_{1}}-{{\gamma }_{1}}-{{\gamma }_{2}} \right), \\
\end{array}
\end{equation*}
or
\begin{equation*}
\begin{array}{l}
     I\left( x,y \right)= \displaystyle \frac{1}{\Gamma \left( {{\gamma }_{1}} \right)\Gamma \left( {{\gamma }_{2}} \right)\Gamma \left( {{\delta }_{1}}-{{\gamma }_{1}}-{{\gamma }_{2}} \right)}\sum\limits_{m,n=0}^{\infty }{{}}\frac{{{\left( {{\gamma }_{3}} \right)}_{{{\beta }_{2}}n}}}{\Gamma \left( {{\delta }_{2}}+{{\alpha }_{4}}m \right)\Gamma \left( {{\delta }_{3}}+{{\beta }_{4}}n \right)}{{x}^{m}}{{y}^{n}} \\
 \times \displaystyle \frac{\Gamma \left( {{\gamma }_{2}}+{{\alpha }_{2}}m \right)\Gamma \left( {{\delta }_{1}}-{{\gamma }_{2}}+\left( {{\alpha }_{3}}-{{\alpha }_{2}} \right)m+{{\beta }_{1}}n \right)}{\Gamma \left( {{\delta }_{1}}+{{\alpha }_{3}}m+{{\beta }_{3}}n \right)}\\
\times \displaystyle \frac{\Gamma \left( {{\gamma }_{1}}+\left( {{\alpha }_{3}}-{{\alpha }_{2}} \right)m+{{\beta }_{1}}n \right)\Gamma \left( {{\delta }_{1}}-{{\gamma }_{1}}-{{\gamma }_{2}} \right)}{\Gamma \left( {{\delta }_{1}}-{{\gamma }_{2}}+\left( {{\alpha }_{3}}-{{\alpha }_{2}} \right)m+{{\beta }_{1}}n \right)}= \\
 = \displaystyle \frac{1}{\Gamma \left( {{\gamma }_{1}} \right)\Gamma \left( {{\gamma }_{2}} \right)}\sum\limits_{m,n=0}^{\infty }{{}}\frac{\Gamma \left( {{\gamma }_{1}}+\left( {{\alpha }_{3}}-{{\alpha }_{2}} \right)m+{{\beta }_{1}}n \right)\Gamma \left( {{\gamma }_{2}}+{{\alpha }_{2}}m \right){{\left( {{\gamma }_{3}} \right)}_{{{\beta }_{2}}n}}}{\Gamma \left( {{\delta }_{1}}+{{\alpha }_{3}}m+{{\beta }_{3}}n \right)}\\
 \times \displaystyle \frac{{{x}^{m}}}{\Gamma \left( {{\delta }_{2}}+{{\alpha }_{4}}m \right)}\frac{{{y}^{n}}}{\Gamma \left( {{\delta }_{3}}+{{\beta }_{4}}n \right)}
 ={{D}_{1}}\left( \begin{matrix}
   {{\gamma }_{1}},{{\alpha }_{3}}-{{\alpha }_{2}},{{\beta }_{1}};{{\gamma }_{2}},{{\alpha }_{2}};{{\gamma }_{3}},{{\beta }_{2}};  \\
   {{\delta }_{1}},{{\alpha }_{3}},{{\beta }_{1}};{{\delta }_{2}},{{\alpha }_{4}};{{\delta }_{3}},{{\beta }_{4}};  \\
\end{matrix}\left| \begin{matrix}
   x  \\
   y  \\
\end{matrix} \right. \right). \\
\end{array}
\end{equation*}
This completes the proof of the integral representation (\ref{eq 21}).
\end{proof}

\section{Laplace transforms }

Let $L_1$ and $L_2$ denote the one-dimensional and two-dimensional Laplace transforms:
\begin{equation*}
\begin{array}{l}
     {{L}_{1}}\left\{ f\left( t \right);p \right\}=  \displaystyle\int\limits_{0}^{\infty }{{}}f\left( t \right){{e}^{-pt}}dt,\,\,\,\operatorname{Re}\,p>0,
\end{array}
\end{equation*}
\begin{equation*}
\begin{array}{l}
    {{L}_{2}}\left\{ f\left( {{t}_{1}},{{t}_{2}} \right);p,q \right\}=  \displaystyle\int\limits_{0}^{\infty }{\int\limits_{0}^{\infty }{{}}}f\left( {{t}_{1}},{{t}_{2}} \right){{e}^{-{{t}_{1}}p-{{t}_{2}}q}}d{{t}_{1}}d{{t}_{2}},\,\,\,\operatorname{Re}\,p>0,\,\,\operatorname{Re}\,q>0.
\end{array}
\end{equation*}
The following Laplace transformations are valid
\begin{equation}\label{eq 22}
\begin{array}{l}
    {{L}_{1}}\left\{ {{t}^{{{\delta }_{1}}-1}}{{D}_{1}}\left( \begin{matrix}
   {{\gamma }_{1}},{{\alpha }_{1}},{{\beta }_{1}};{{\gamma }_{2}},{{\alpha }_{2}};{{\gamma }_{3}},{{\beta }_{2}};  \\
   {{\delta }_{1}},{{\alpha }_{3}},{{\beta }_{3}};{{\delta }_{2}},{{\alpha }_{4}};{{\delta }_{3}},{{\beta }_{4}};  \\
\end{matrix}\left| \begin{matrix}
   x{{t}^{{{\alpha }_{3}}}}  \\
   y{{t}^{{{\beta }_{3}}}}  \\
\end{matrix} \right. \right):p \right\} \\
= \displaystyle \frac{1}{{{p}^{{{\delta }_{1}}}}}{{E}_{11}}\left( \begin{matrix}
   {{\gamma }_{1}},{{\alpha }_{1}},{{\beta }_{1}};{{\gamma }_{2}},{{\alpha }_{2}};{{\gamma }_{3}},{{\beta }_{2}};  \\
   {{\delta }_{2}},{{\alpha }_{4}};{{\delta }_{3}},{{\beta }_{4}};  \\
\end{matrix}\left| \begin{matrix}
   \displaystyle \frac{x}{{{p}^{{{\alpha }_{3}}}}}  \\
   \displaystyle \frac{y}{{{p}^{{{\beta }_{3}}}}}  \\
\end{matrix} \right. \right), \operatorname{Re}p>0, \\
\end{array}
\end{equation}
\begin{equation}\label{eq 23}
\begin{array}{l}
   {{L}_{1}}\left\{ {{t}^{\rho -1}}{{E}_{1}}\left( \begin{matrix}
   {{\gamma }_{1}},{{\alpha }_{1}};{{\gamma }_{2}},{{\beta }_{1}};  \\
   {{\delta }_{1}},{{\alpha }_{2}},{{\beta }_{2}};{{\delta }_{2}},{{\alpha }_{3}};{{\delta }_{3}},{{\beta }_{3}};  \\
\end{matrix}\left| \begin{matrix}
   x{{t}^{{{\mu }_{1}}}}  \\
   y{{t}^{{{\mu }_{2}}}}  \\
\end{matrix} \right. \right);p \right\} \\
 = \displaystyle \frac{\Gamma \left( \rho  \right)}{{{p}^{\rho }}}{{D}_{1}}\left( \begin{matrix}
   \rho ,{{\mu }_{1}},{{\mu }_{2}};{{\gamma }_{1}},{{\alpha }_{1}};{{\gamma }_{2}},{{\beta }_{1}};  \\
   {{\delta }_{1}},{{\alpha }_{2}},{{\beta }_{2}};{{\delta }_{2}},{{\alpha }_{3}};{{\delta }_{3}},{{\beta }_{3}};  \\
\end{matrix}\left| \begin{matrix}
   \displaystyle \frac{x}{{{p}^{{{\mu }_{1}}}}}  \\
   \displaystyle \frac{y}{{{p}^{{{\mu }_{2}}}}}  \\
\end{matrix} \right. \right),\\
\operatorname{Re}\rho >0,\,\operatorname{Re}{{\mu }_{1}}>0,\,\operatorname{Re}{{\mu }_{2}}>0, \\
\end{array}
\end{equation}
\begin{equation}\label{eq 24}
\begin{array}{l}
  {{L}_{2}}\left\{ t_{1}^{{{\rho }_{1}}-1}t_{2}^{{{\rho }_{2}}-1}{{E}_{8}}\left( \begin{matrix}
   {{\gamma }_{1}},{{\alpha }_{1}},{{\beta }_{1}};  \\
   {{\delta }_{1}},{{\alpha }_{2}},{{\beta }_{2}};{{\delta }_{2}},{{\alpha }_{3}};{{\delta }_{3}},{{\beta }_{3}};  \\
\end{matrix}\left| \begin{matrix}
   xt_{1}^{{{\mu }_{1}}}  \\
   yt_{2}^{{{\mu }_{2}}}  \\
\end{matrix} \right. \right);p,q \right\}= \\
 = \displaystyle \frac{\Gamma \left( {{\rho }_{1}} \right)\Gamma \left( {{\rho }_{2}} \right)}{{{p}^{{{\rho }_{1}}}}{{q}^{{{\rho }_{2}}}}}{{D}_{1}}\left( \begin{matrix}
   {{\gamma }_{1}},{{\alpha }_{1}},{{\beta }_{1}};{{\rho }_{1}},{{\mu }_{1}};{{\rho }_{2}},{{\mu }_{2}};  \\
   {{\delta }_{1}},{{\alpha }_{2}},{{\beta }_{2}};{{\delta }_{2}},{{\alpha }_{3}};{{\delta }_{3}},{{\beta }_{3}};  \\
\end{matrix}\left| \begin{matrix}
   x  \\
   y  \\
\end{matrix} \right. \right),\,\,\,\operatorname{Re}{{\rho }_{1}}>0,\,\,\,\operatorname{Re}{{\rho }_{2}}>0, \\
\end{array}
\end{equation}
\begin{equation}\label{eq 25}
\begin{array}{l}
   {{L}_{2}}\left\{ {{x}^{{{\delta }_{1}}-1}}{{y}^{{{\delta }_{2}}-1}}{{D}_{2}}\left( \begin{matrix}
   {{\gamma }_{1}},{{\alpha }_{1}},{{\beta }_{1}};{{\gamma }_{2}},{{\alpha }_{2}};{{\gamma }_{3}},{{\beta }_{2}};  \\
   {{\delta }_{1}},{{\alpha }_{3}};{{\delta }_{2}},{{\beta }_{3}};{{\delta }_{3}},{{\alpha }_{4}};{{\delta }_{4}},{{\beta }_{4}};  \\
\end{matrix}\left| \begin{matrix}
   {{x}^{{{\alpha }_{3}}}}  \\
   {{y}^{{{\beta }_{3}}}}  \\
\end{matrix} \right. \right);p,q \right\} \\
 =\displaystyle \frac{1}{{{p}^{{{\delta }_{1}}}}{{q}^{{{\delta }_{2}}}}}{{E}_{11}}\left( \begin{matrix}
   {{\gamma }_{1}},{{\alpha }_{1}},{{\beta }_{1}};{{\gamma }_{2}},{{\alpha }_{2}};{{\gamma }_{3}},{{\beta }_{2}};  \\
   {{\delta }_{3}},{{\alpha }_{4}};{{\delta }_{4}},{{\beta }_{4}};  \\
\end{matrix}\left| \begin{matrix}
  \displaystyle \frac{1}{{{p}^{{{\alpha }_{3}}}}}  \\
  \displaystyle \frac{1}{{{q}^{{{\beta}_{3}}}}}  \\
\end{matrix} \right. \right), \operatorname{Re}p>0,\,\,\operatorname{Re}q>0, \\
\end{array}
\end{equation}

Equalities (\ref{eq 22})-(\ref{eq 25}) can be verified by direct calculations.

\section{System of partial differential equations}

{\bf Theorem. }
        \label{thefirstone}
       Let $\theta \equiv x\left( \partial /\partial x \right)$, $\phi \equiv y\left( \partial /\partial y \right)$, then for ${{\gamma }_{1}},{{\gamma }_{2}},{{\gamma }_{3}},{{\delta }_{1}},{{\delta }_{2}},$ ${{\delta }_{3}},x,y\in \mathbb{C}$, ${{\alpha }_{j}},{{\beta }_{j}}\in \mathbb{N},\,\,\,\left( j=1,2,3,4 \right)$   the two variable Mittag-Leffler-type function $D_1$, defined by (\ref{eq 2.1}), satisfies the following partial differential equations:
\begin{equation}\label{eq 26}
\begin{array}{l}
       \prod\limits_{i=1}^{{{\alpha }_{3}}}{{}}\left( {{\delta }_{1}}+{{\alpha }_{3}}-i+{{\alpha }_{3}}\theta +{{\beta }_{3}}~\phi  \right)\prod\limits_{i=1}^{{{\alpha }_{4}}}{{}}\left( {{\delta }_{2}}+{{\alpha }_{4}}-i+{{\alpha }_{4}}\theta  \right){{x}^{-1}} \\
 -\prod\limits_{i=1}^{{{\alpha }_{1}}}{{}}\left( {{\gamma }_{1}}+{{\alpha }_{1}}-i+{{\alpha }_{1}}\theta +{{\beta }_{1}}\phi  \right)\prod\limits_{i=1}^{{{\alpha }_{2}}}\left( {{\gamma }_{2}}+{{\alpha }_{2}}-i+{{\alpha }_{2}}\theta  \right){{D}_{1}}(x,y)=0,\\
  \prod\limits_{i=1}^{{{\beta }_{3}}}{{}}\left( {{\delta }_{1}}+{{\beta }_{3}}-i+{{\alpha }_{3}}\theta +{{\beta }_{3}}~\phi  \right)\prod\limits_{i=1}^{{{\beta }_{4}}}\left( {{\delta }_{2}}+{{\beta }_{4}}-i+{{\beta }_{4}}~\phi  \right){{y}^{-1}} \\
  -\prod\limits_{i=1}^{{{\beta }_{1}}}{{}}\left( {{\gamma }_{1}}+{{\beta }_{1}}-i+{{\alpha }_{1}}\theta +{{\beta }_{1}}~\phi  \right)\prod\limits_{i=1}^{{{\beta }_{2}}}\left( {{\gamma }_{2}}+{{\beta }_{2}}-i+{{\beta }_{2}}~\phi  \right) {{D}_{1}}(x,y)=0.
\end{array}
\end{equation}

\begin{proof}
Let us show the validity of the first equation in (\ref{eq 26}). The right side of the first equation (\ref{eq 26}) will be denoted by $J(x,y)$. Substituting the function $D_1$  into the first equation (\ref{eq 26}) and considering the equalities $\theta \left( {{x}^{m}} \right)=m{{x}^{m}},\,\,\phi \left( {{y}^{n}} \right)=n{{y}^{n}},$ we obtain
\begin{equation*}
\begin{array}{l}
        \prod\limits_{i=1}^{{{\alpha }_{4}}}{{}}\left( {{\delta }_{2}}+{{\alpha }_{4}}-i+{{\alpha }_{4}}\theta  \right){{x}^{-1}}{{D}_{1}}(x,y) \\
  = \displaystyle \sum\limits_{m=1}^{\infty }{\sum\limits_{n=0}^{\infty }{\frac{{{\left( {{\gamma }_{1}} \right)}_{{{\alpha }_{1}}m+{{\beta }_{1}}n}}{{\left( {{\gamma }_{2}} \right)}_{{{\alpha }_{2}}m}}{{\left( {{\gamma }_{3}} \right)}_{{{\beta }_{2}}n}}}{\Gamma \left( {{\delta }_{1}}+{{\alpha }_{3}}m+{{\beta }_{3}}n \right)}\frac{{{x}^{m-1}}}{\Gamma \left( {{\delta }_{2}}+{{\alpha }_{4}}\left( m-1 \right) \right)}\frac{{{y}^{n}}}{\Gamma \left( {{\delta }_{3}}+{{\beta }_{4}}n \right)}}}, \\
 \prod\limits_{i=1}^{{{\alpha }_{2}}}{{}}\left( {{\gamma }_{2}}+{{\alpha }_{2}}-i+{{\alpha }_{2}}\theta  \right){{D}_{1}}(x,y) \\
 =\displaystyle \sum\limits_{m,n=0}^{\infty }{\frac{{{\left( {{\gamma }_{1}} \right)}_{{{\alpha }_{1}}m+{{\beta }_{1}}n}}{{\left( {{\gamma }_{2}} \right)}_{{{\alpha }_{2}}\left( m+1 \right)}}{{\left( {{\gamma }_{3}} \right)}_{{{\beta }_{2}}n}}}{\Gamma \left( {{\delta }_{1}}+{{\alpha }_{3}}m+{{\beta }_{3}}n \right)}\frac{{{x}^{m}}}{\Gamma \left( {{\delta }_{2}}+{{\alpha }_{4}}m \right)}\frac{{{y}^{n}}}{\Gamma \left( {{\delta }_{3}}+{{\beta }_{4}}n \right)}}. \\
\end{array}
\end{equation*}
Hence we get
\begin{equation}\label{eq 27}
\begin{array}{l}
  \prod\limits_{i=1}^{{{\alpha }_{3}}}{{}}\left( {{\delta }_{1}}+{{\alpha }_{3}}-i+{{\alpha }_{3}}\theta +{{\beta }_{3}}~\phi  \right)\prod\limits_{i=1}^{{{\alpha }_{4}}}{{}}\left( {{\delta }_{2}}+{{\alpha }_{4}}-i+{{\alpha }_{4}}\theta  \right){{x}^{-1}}{{D}_{1}}(x,y) \\
 = \displaystyle \sum\limits_{m=1}^{\infty }{\sum\limits_{n=0}^{\infty }{\frac{{{\left( {{\gamma }_{1}} \right)}_{{{\alpha }_{1}}m+{{\beta }_{1}}n}}{{\left( {{\gamma }_{2}} \right)}_{{{\alpha }_{2}}m}}{{\left( {{\gamma }_{3}} \right)}_{{{\beta }_{2}}n}}}{\Gamma \left( {{\delta }_{1}}+{{\alpha }_{3}}\left( m-1 \right)+{{\beta }_{3}}n \right)}\frac{{{x}^{m-1}}}{\Gamma \left( {{\delta }_{2}}+{{\alpha }_{4}}\left( m-1 \right) \right)}\frac{{{y}^{n}}}{\Gamma \left( {{\delta }_{3}}+{{\beta }_{4}}n \right)},}} \\
\end{array}
\end{equation}
\begin{equation}\label{eq 28}
\begin{array}{l}
     \prod\limits_{i=1}^{{{\alpha }_{1}}}{{}}\left( {{\gamma }_{1}}+{{\alpha }_{1}}-i+{{\alpha }_{1}}\theta +{{\beta }_{1}}\phi  \right)\prod\limits_{i=1}^{{{\alpha }_{2}}}{{}}\left( {{\gamma }_{2}}+{{\alpha }_{2}}-i+{{\alpha }_{2}}\theta  \right){{D}_{1}}(x,y) \\
 = \displaystyle \sum\limits_{m,n=0}^{\infty }{\frac{{{\left( {{\gamma }_{1}} \right)}_{{{\alpha }_{1}}\left( m+1 \right)+{{\beta }_{1}}n}}{{\left( {{\gamma }_{2}} \right)}_{{{\alpha }_{2}}\left( m+1 \right)}}{{\left( {{\gamma }_{3}} \right)}_{{{\beta }_{2}}n}}}{\Gamma \left( {{\delta }_{1}}+{{\alpha }_{3}}m+{{\beta }_{3}}n \right)}\frac{{{x}^{m}}}{\Gamma \left( {{\delta }_{2}}+{{\alpha }_{4}}m \right)}\frac{{{y}^{n}}}{\Gamma \left( {{\delta }_{3}}+{{\beta }_{4}}n \right)}}. \\
\end{array}
\end{equation}
Substituting (\ref{eq 27}) - (\ref{eq 28}) into the first equation of system (\ref{eq 26}), we determine
\begin{equation*}
\begin{array}{l}
J(x,y)= \displaystyle \sum\limits_{m=1}^{\infty }{\sum\limits_{n=0}^{\infty }{\frac{{{\left( {{\gamma }_{1}} \right)}_{{{\alpha }_{1}}m+{{\beta }_{1}}n}}{{\left( {{\gamma }_{2}} \right)}_{{{\alpha }_{2}}m}}{{\left( {{\gamma }_{3}} \right)}_{{{\beta }_{2}}n}}}{\Gamma \left( {{\delta }_{1}}+{{\alpha }_{3}}\left( m-1 \right)+{{\beta }_{3}}n \right)}\frac{{{x}^{m-1}}}{\Gamma \left( {{\delta }_{2}}+{{\alpha }_{4}}\left( m-1 \right) \right)}\frac{{{y}^{n}}}{\Gamma \left( {{\delta }_{3}}+{{\beta }_{4}}n \right)}-}} \\
 - \displaystyle \sum\limits_{m,n=0}^{\infty }{\frac{{{\left( {{\gamma }_{1}} \right)}_{{{\alpha }_{1}}\left( m+1 \right)+{{\beta }_{1}}n}}{{\left( {{\gamma }_{2}} \right)}_{{{\alpha }_{2}}\left( m+1 \right)}}{{\left( {{\gamma }_{3}} \right)}_{{{\beta }_{2}}n}}}{\Gamma \left( {{\delta }_{1}}+{{\alpha }_{3}}m+{{\beta }_{3}}n \right)}\frac{{{x}^{m}}}{\Gamma \left( {{\delta }_{2}}+{{\alpha }_{4}}m \right)}\frac{{{y}^{n}}}{\Gamma \left( {{\delta }_{3}}+{{\beta }_{4}}n \right)}}. \\
\end{array}
\end{equation*}
In the first term, changing the summation index $m$  to $m+1$, we find that $J(x,y)=0$.
The validity of the second equation of  (\ref{eq 26}) can be proved similarly. 
\end{proof}

\section{Conclusion} 
We have introduced the series of two-variable Mittag-Leffler-type functions and studied certain properties of these functions. Namely, we have determined the region of the convergence, Euler-type integral representation, and one and two-dimensional Laplace transform and determined the system of partial differential equations linked with these functions. We note that we did these for the function $D_1(x,y)$, but results can be easily obtained for other introduced functions. We believe that the obtained result will be applied in the near future since similar functions $E_1(x,y)$ and $E_2(x,y)$ already found their applications.  

\section{Conflict of Interests}
This work does not have any conflicts of interest.

\end{document}